\documentclass[a4paper,10pt]{article}
\usepackage[cp1251]{inputenc}
\usepackage[english]{babel}

\usepackage{amssymb}\usepackage{amsmath}
\usepackage{epsfig}
\usepackage{psfig}
\def\bl{\rule[-1mm]{2.4mm}{2.4mm}}

\def\be{\begin{equation}}
\def\ee{\end{equation}}

\newtheorem{thrm}{\bf Theorem}
\newtheorem{lmm}{\bf Lemma}
\newtheorem{dfn}{\bf Definition}
\newtheorem{rmk}{\bf Remark}

\begin{document}

\title {Conformal mapping of rectangular heptagons}
\author{\copyright 2011 ~~~~A.B.Bogatyrev
\thanks{Supported by RFBR grants 10-01-00407 and RAS Program
"Modern problems of theoretical mathematics"}}
\date{}
\maketitle

\section{Introduction}
There exists an impressive list of numerical methods for conformal mapping
of polygons \cite{TD}. A simple observation allows one to extend this list.
Once the angles of the polygon are rational multiples of
$\pi$, the Christoffel-Schwartz (CS) integral which maps the upper half plane
 $\mathbb{H}$ to the polygon is an abelian integral on a compact Riemann surface. The full power of the function theory on Riemann surfaces may be applied now to attack the evaluation of the CS integral as well as its auxiliary parameters.

In this paper we consider the simplest case beyond the elliptic
one (4 right angles in the polygon) described in \cite{BHY} and elaborate this
approach to the case of a simply connected polygon with six right and one zero
angle. Two of the right angles have to be exterior ones (i.e.
equal to $3\pi/2$) and the vertex with zero angle is infinitely
distant. We give the representation for the  both mappings,
heptagon to the half plane and back as certain explicit
expressions in terms of  genus two Riemann theta functions. For
the latter there exist a robust and effective method of
computation \cite{DHB} with controllable accuracy. Therefore we
can guarantee the machine accuracy for  the conformal mapping uniformly in the polygon/halfplane.

As usual, there are several auxiliary parameters of the mapping which are determined by the geometrical dimensions
of the polygon. A good portion of those determining equations are linear in  this approach.  Say, if we map the L-shaped rectangular hexagon, one has to solve just one nonlinear equation to determine the mapping.

This method may be used as a reference for testing the numerical conformal mappings. It works equally well for the non-convex polygons, in the presence of  spikes of the boundary, narrow isthma, boundary elements of different scales, etc.

\section{Spaces of heptagons}
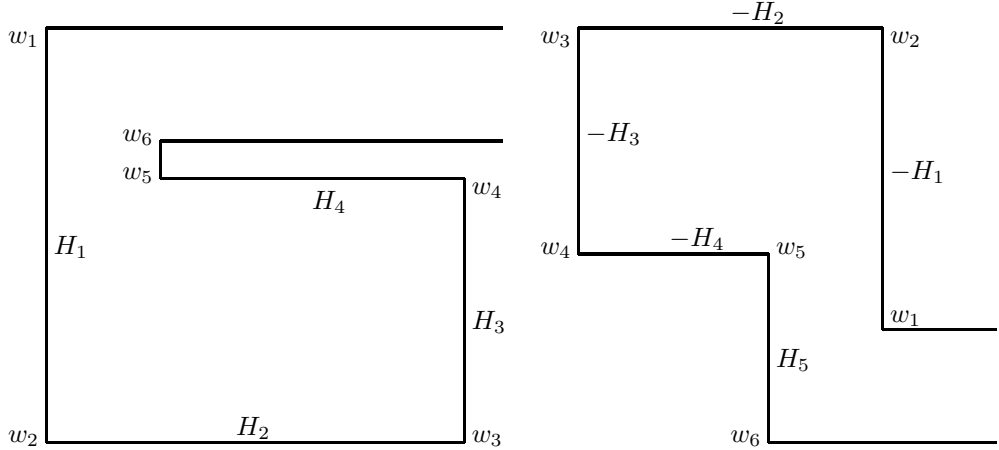
\begin{figure}[h]
\begin{picture}(170,65)
\thicklines
\put(5,60){\line(1,0){60}}
\put(5,60){\line(0,-1){55}}
\put(5,5){\line(1,0){55}}
\put(60,5){\line(0,1){35}}
\put(60,40){\line(-1,0){40}}
\put(20,40){\line(0,1){5}}
\put(20,45){\line(1,0){45}}
\put(0,58){$w_1$}
\put(6,30){$H_1$}
\put(0,5){$w_2$}
\put(30,6){$H_2$}
\put(61,5){$w_3$}
\put(61,20){$H_3$}
\put(61,38){$w_4$}
\put(40,36){$H_4$}
\put(15,40){$w_5$}
\put(15,45){$w_6$}

\put(100,5){\line(1,0){30}}
\put(100,30){\line(0,-1){25}}
\put(75,30){\line(1,0){25}}
\put(75,60){\line(0,-1){30}}
\put(115,60){\line(-1,0){40}}
\put(115,20){\line(0,1){40}}
\put(115,20){\line(1,0){15}}
\put(116,21){$w_1$}
\put(116,40){$-H_1$}
\put(116,58){$w_2$}
\put(95,61){$-H_2$}
\put(70,58){$w_3$}
\put(76,45){$-H_3$}
\put(70,30){$w_4$}
\put(87,31){$-H_4$}
\put(101,30){$w_5$}
\put(101,15){$H_5$}
\put(95,5){$w_6$}

\end{picture}
\caption{\small Heptagons from spaces ${\cal P}_{56}$ (left) and  ${\cal P}_{15}$ (right) }
\label{Heptagons}
\end{figure}

We consider rectangular heptagon with half-infinite  width $\pi$ "channel" oriented to the east as it is shown in the Fig. \ref{Heptagons}. Its sides are either vertical or horizontal. The  vertex of the polygon at infinity is denoted $w_0$, others  are enumerated in increasing order with respect to the natural orientation of the boundary (counterclockwise). Two vertexes with the angles $3\pi/2$  are given special names $w_\alpha$ and $w_\beta$, $1\le\alpha<\beta\le6$.

 The space ${\cal P}_{\alpha\beta}$ of the heptagons
 with the given indexes $\alpha,~\beta$ is parametrized by the lengths of the sides, to which we ascribe signs for technical reasons:
\be
\label{HeptDim}
i^sH_s:=w_s-w_{s+1}\qquad s=1,2,\dots,5;\\
\ee
One obvious restriction on the real values $H_1,\dots,H_5$ is the following
\be
H_1-H_3+H_5=\pi \quad (=:Im(w_1-w_6)).
\label{sumH}
\ee
 The sign of $H_s$ is negative iff $\alpha\le s<\beta$:
\be
\label{SignRule}
(s+\frac12-\alpha)(s+\frac12-\beta)~H_s>0.
\ee
The boundary of the heptagon has no self-intersections which means that the dimensions obey the following additional inequalities:
\be
\label{NoIntersection}
\begin{array}{c|cc}
(\alpha,\beta)&{\rm Restriction}&\\
\hline\\
(1,2)&-H_2+H_4>0&\\
(1,5)&-H_2+H_4>0&when ~H_1-H_3\le0\\
(2,3)&-H_3+H_5>0&\\
(2,6)&-H_2+H_4<0&when -H_3+H_5\le0\\
(4,5)&H_1-H_3>0&\\
(5,6)&-H_2+H_4<0&\\
\end{array}
\ee

\begin{lmm}\label{H1H5}
The heptagons with fixed indexes $\alpha$ and $\beta$ make up a
connected space ${\cal P}_{\alpha\beta}$ of real dimension 4 with the global
coordinates $H_1,\dots,H_5$ subjected to restriction
(\ref{sumH}), the sign rule (\ref{SignRule}) and inequalities
(\ref{NoIntersection}). ~~~\bl
\end{lmm}

Of course, the point from the space ${\cal P}_{\alpha,\beta}$ defines a heptagon up to translations in the plane only. Those translations may be eliminated when necessary by the normalization e.g. $w_1:=i\pi$. The reflection in the real axis induces the mapping ${\cal P}_{\alpha,\beta}\to$ ${\cal P}_{7-\beta,7-\alpha}$
which in the above coordinate system appears as $(H_1,H_2,\dots,H_5)\to$ $(H_5,H_4,\dots,H_1)$.

\section{Hyperelliptic curves with six real branch points}
The conformal mapping of the upper half plane to any heptagon
from the space ${\cal P}_{\alpha\beta}$ may be represented by the Christoffel-Schwarz integral. This integral lives on a hyperelliptic curve
with six real branchpoints. In this section we briefly remind several facts about such curves.

\subsection{Algebraic model}
The double cover of the sphere with six real branch points
$x_1<x_2<\dots<x_5<x_6$
 is a compact genus two Riemann surface $X$ with the equation
 (of its affine part):
\be
\label{X}
y^2=\prod\limits_{s=1}^6(x-x_s),
\qquad  (x,y)\in\mathbb{C}^2.
\ee
 This curve admits a \emph{conformal involution} $J(x,y)=(x,-y)$ with six stationary points $p_s:=(x_s,0)$ and an \emph{anticonformal involution} (\emph{reflection})
$\bar{J}(x,y)=(\bar{x},\bar{y})$. The stationary points set of the latter has three components known as \emph{real ovals} of the curve. Each real oval is an embedded circle \cite{Nat} and doubly covers exactly one of the segments  $[x_2,x_3]$, $[x_4,x_5]$, $[x_6,x_1]\ni\infty$ of the extended real line $\hat{\mathbb R}:={\mathbb R}\cup\infty$. We denote those ovals as \emph{first, second and third} respectively. The  lift of the complimentary set of intervals
$[x_1,x_2]$, $[x_3,x_4]$, $[x_5,x_6]$ to the surface (\ref{X})
 gives us three \emph{coreal ovals} which make up the set of points fixed by another anticonformal involution ${\bar J}J=J{\bar J}$.

\subsection{Homologies, differentials, periods}\label{SectBasisCycles}
\begin{figure}
\psfig{file=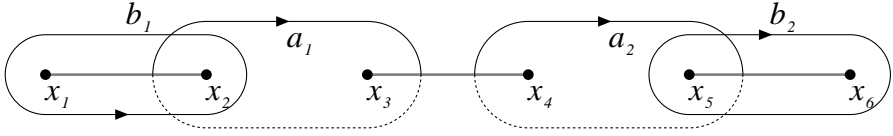}
\caption{\small Canonical basis in homologies of the curve $X$}
\label{BasisHomologies}
\end{figure}

We fix a special basis in the 1-homology space of the curve $X$ intrinsically related to the latter. The first and second real ovals give us two 1-cycles,  $a_1$ and $a_2$ respectively. Both cycles are oriented (up to simultaneous change of sign) as the boundary of a pair of pants obtained by removing real ovals from the surface. Two remaining cycles $b_1$ and $b_2$ are coreal ovals of the curve oriented so that the intersection matrix takes the canonical form --
see Fig. \ref{BasisHomologies}.

The reflection of the surface acts on the introduced basis as follows
\be
\label{abreflect}
\begin{array}{c}
\bar{J}a_s=a_s,
\quad
\bar{J}b_s=-b_s,
\end{array}
\qquad s=1,2.
\ee
Holomorphic differentials on the curve $X$ take the form
\be
\label{diffRep}
du_*=(C_{1*}x+C_{2*})y^{-1}dx,
\ee
with constant values $C_{1*}$, $C_{2*}$.
 The basis of differentials dual to the basis of cycles
\be
\int_{a_s}du_j:=\delta_{sj};
\qquad s,j =1,2,
\ee
determines Riemann period matrix $\Pi$ with the elements
\be
\label{periods}
\Pi_{sj}:=\int_{b_s}du_j;
\qquad s,j =1,2.
\ee
It is a classical result that $\Pi$ is symmetric and has positive definite imaginary part \cite{FK}.

From the symmetry properties (\ref{abreflect}) of the chosen basic cycles it readily follows that:
\begin{itemize}
 \item Normalized differentials are real ones, i.e. $\bar{J}du_s=\overline{du_s}$, in other words the coefficients $C_*$ in the representation (\ref{diffRep}) are real.
  \item Period matrix is purely imaginary,  therefore we can introduce the symmetric and positive definite real  matrix  $\Omega:=Im(\Pi)$,
 \item Zeroes of the differential $du_2$ (resp. $du_1$) lie on the first (resp. second) real oval.
\end{itemize}

\subsection{ Jacobian and Abel-Jacobi map }
\begin{dfn}
Given a Riemann period matrix $\Pi$, we define the full rank lattice
\be
\label{Lattice}
L(\Pi)=\Pi\mathbb{Z}^2+ \mathbb{Z}^2 = \int_{H_1(X,\mathbb{Z})}du,
\qquad du:=(du_1,du_2)^t,
\ee
in $\mathbb{C}^2$ and the 4-torus $Jac(X):=\mathbb{C}^2/L(\Pi)$ known as a Jacobian of the curve $X$.
\end{dfn}
This definition depends on the choice of the basis in the lattice $H_1(X,\mathbb{Z})$, other choices bring us to isomorphic tori.

It is convenient to represent the points  $u\in\mathbb{C}^2$ as a theta characteristic $[\epsilon,\epsilon']$, i.e. a couple of real 2-vectors (columns) $\epsilon, \epsilon'$:
\be
u=\frac12(\Pi\epsilon+\epsilon').
\ee
The points of Jacobian  $Jac(X)$ in this notation correspond to
two vectors with real entries modulo 2. Second order points of Jacobian are
represented as $2\times 2$ matrices with ${\mathbb Z}_2$ entries.
Please note that we use a nonstandard notation of theta characteristic as two column vectors written one after another. Usually the transposed matrix is used.

\begin{dfn}
Abel-Jacobi (briefly: AJ) map is a correctly defined mapping from the surface $X$ to its Jacobian.
\be
\label{AJmap}
u(p):=\int_{p_1}^p du~~ mod~L(\Pi),
\qquad p_1:=(x_1,0); \quad du:=(du_1,du_2)^t,
\ee
\end{dfn}

From Riemann-Roch formula it easily follows \cite{FK} that Abel-Jacobi map is a holomorphic embedding of the curve into its Jacobian. In Sect. \ref{theta} we give an explicit equation for the image of the genus two curve in its Jacobian. Let us meanwhile compute the images of the branching points $p_s=(x_s,0)$, $s=1,\dots,6$ of the curve $X$:

\centerline{
\begin{tabular}{c|c|c}
$p$ & $u(p)~ mod~ L(\Pi)$ & $[\epsilon, \epsilon'](u(p))$\\
\hline
$p_2$ & $\Pi^1/2$ & $\tiny
\left[\begin{array}{c} 10\\00\end{array}\right]$\\
$p_3$&$(\Pi^1+E^1)/2$&$\tiny
\left[\begin{array}{c} 11\\00\end{array}\right]$\\
$p_4$&$(\Pi^2+E^1)/2$&$\tiny
\left[\begin{array}{c} 01\\10\end{array}\right]$\\
$p_5$&$(\Pi^2+E^1+E^2)/2$&$\tiny
\left[\begin{array}{c}01\\11\end{array}\right]$\\
$p_6$&$(E^1+E^2)/2$&$\tiny
\left[\begin{array}{c}01\\01\end{array}\right]$\\
\label{AJPj}
\end{tabular}}

where $\Pi^s$ and $E^s$ are the $s$-th columns of the period and identity matrix respectively. One can notice that vector $\epsilon(u(p))$ is constant along the real ovals and $\epsilon'(u(p))$ is constant along the coreal ovals.

\subsection{Tiling the Jacobian}
 Let us consider  sixteen  disjoint blocks (tiles) in the Jacobian of a curve filled  by the points with theta characteristics in the sets
$$
\left[
\begin{array}{c}
\pm I ~~\pm I\\
\pm I ~~\pm I
\end{array}
\right],
\qquad I:=(0,1),
$$
distinguished by the possible choices of four signs.
Jacobian itself  is the closure of the  union of those tiles.

 Chosen an orientation on a real oval of the curve, one can distinguish a component in the set $x^{-1}\mathbb{H}\subset X$ which lies to the left of this oval. We shall denote this disc as $\mathbb{H}^+$.
 The surface $X$ with all real and coreal ovals removed is a disjoint union of four open 2-discs $\mathbb{H}^+$, $J\mathbb{H}^+$, $\bar{J}\mathbb{H}^+$, $\bar{J}J\mathbb{H}^+$.  It turns out that AJ map sends each of those discs to a certain block of the Jacobian and each (co)real oval -- to a certain 2-torus. This helps us to discriminate points $p\in X$ with the same value of the projection $x(p)$

 \begin{thrm}
 Let the disc $\mathbb{H}^+$ be chosen in accordance with the orientation of
 (any of) the $a$-cycles. Then the above four discs are mapped to the following four blocks of the Jacobian:
\end{thrm}
\centerline{
\begin{tabular}{c|c|c|c|c}
\label{RangeDisc}
$p\in$ & $\mathbb{H}^+$ & $J\mathbb{H}^+$ & $\bar{J}\mathbb{H}^+$& $\bar{J}J\mathbb{H}^+$\\
\hline
$[\epsilon,\epsilon'](u(p))\in$ &
$\tiny\left[
\begin{array}{cc}
-I~&I\\
-I~&I\\
\end{array}
\right]$
&
$\tiny\left[
\begin{array}{cc}
I~&-I\\
I~&-I\\
\end{array}
\right]$
&
$\tiny\left[
\begin{array}{cc}
I~&I\\
I~&I\\
\end{array}
\right]$
&
$\tiny\left[
\begin{array}{cc}
-I~&-I\\
-I~&-I\\
\end{array}
\right]$\\
\end{tabular}}
\vspace{3mm}

{\bf Proof.}
  Symmetries of the normalized abelian differentials with respect to involutions $J$, $\bar{J}$ guarantee the following equalities:
$u(Jp)=-u(p)$; $u(\bar{J}p)=\overline{u(p)}$  since the base point $p_1$ of AJ map is fixed by both involutions. Therefore we may thoroughly investigate the map on the disc $\mathbb{H}^+$ only.

\begin{figure}
\begin{picture}(160,40)
\thicklines
\put(10,5){\vector(0,1){25}}
\put(10,5){\line(1,0){30}}
\put(40,5){\line(0,1){10}}
\put(40,15){\line(-1,0){20}}
\put(20,15){\line(0,1){1}}
\put(20,16){\line(1,0){20}}
\put(40,16){\vector(0,1){14}}
\put(40,30){\line(-1,0){30}}
\put(5,1){$u_1(x_2)$}
\put(35,1){$u_1(x_3)$}
\put(41,13){$u_1(x_4)$}
\put(41,16){$u_1(x_5)$}
\put(35,31){$u_1(x_6)$}
\put(5,31){$u_1(x_1)$}
\put(0,15){$\Omega_{11}/2$}
\put(41,23){$\Omega_{12}/2$}
\put(24,0){$1/2$}
\put(20,21){$u_1(\mathbb{H}^+)$}

\put(90,5){\line(0,1){5}}
\put(90,10){\line(1,0){15}}
\put(105,10){\line(0,1){1}}
\put(105,11){\line(-1,0){15}}
\put(90,11){\vector(0,1){24}}
\put(90,35){\line(1,0){30}}
\put(120,5){\vector(0,1){30}}
\put(90,5){\line(1,0){30}}
\put(85,1){$u_2(x_4)$}
\put(115,1){$u_2(x_5)$}
\put(115,36){$u_2(x_6)$}
\put(85,36){$u_2(x_1)$}
\put(78,11){$u_2(x_2)$}
\put(78,8){$u_2(x_3)$}
\put(121,20){$\Omega_{22}/2$}
\put(80,23){$\Omega_{21}/2$}
\put(105,0){$1/2$}
\put(103,20){$u_2(\mathbb{H}^+)$}

\end{picture}
\caption{\small The images of the upper half plane $u_1(\mathbb{H}^+)$ (left) and  $u_2(\mathbb{H}^+)$ (right)}
\label{URect}
\end{figure}
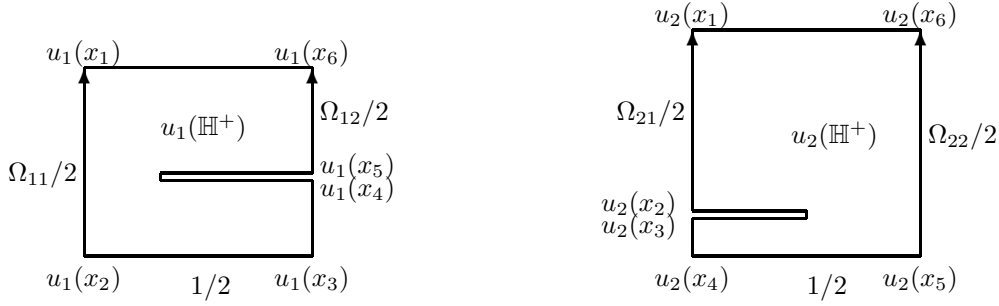

Both components $u_s(p)$ of Abel-Jacobi map are CS integrals and send the disc $\mathbb{H}^+$ to the rectangles with slots shown in the Fig. \ref{URect}. Clearly, the real part of $u_s(p)$ lies in the interval $(0,1/2)$ which may be reformulated as $\epsilon'(u(p))\in (I,I)^t$.

To study the range of 2-vector $\epsilon:=\Omega^{-1}Im~2u(p)$ we introduce new differentials $dv=(dv_1,dv_2)^t:=-i\Omega^{-1}du$. Those are purely imaginary with normalization $\int_{b_j}dv_s=\delta_{js}$. Differential $dv_1$ has zeroes on coreal oval covering the interval $(x_5,x_6)$; $dv_2$ has zeroes which project to the interval $(x_1,x_2)$.
Again, each abelian integral $v_s(p)$ maps the upper half plane to the rectangle with slot similar to that in the left Fig. \ref{URect}, but now the vertex $v_s(x_1)$ lies on the right side of the rectangle, therefore
$\epsilon(u(p)):=-2Re ~v(p)\in (I,I)^t$.
~~~\bl

\section{Two moduli spaces}
We are going to describe the space of CS integrals corresponding to the rectangular polygons from the spaces ${\cal P}_{\alpha\beta}$. This space is an extension of the underlying space of genus 2 curves with real branch points.

\subsection{Genus 2 curves with three real ovals}
Each genus two Riemann surface $X$ is  automatically hyperelliptic, i.e. it admits a conformal involution $J$ that fixes six points. This involution is unique. Riemann surface is said to have a reflection iff it admits an anticonformal involution $\bar{J}$  (same surface $X$ may have several anticonformal involutions). Each component of the set of $\bar{J}$-fixed points  is an embedded circle \cite{Nat} known as a \emph{ real oval} of the reflection.

\begin{dfn} Moduli space ${\cal M}_2\mathbb{R}_3$ is the space of
genus two Riemann surfaces with reflection and three enumerated real ovals. Two surfaces are equivalent iff  there is a conformal $1-1$ mapping between them commuting with the reflections and respecting the marking of real ovals. \end{dfn}

Let us consider the constructive model of an element $X\in {\cal M}_2\mathbb{R}_3$.

Necessarily, the involutions $J$ and $\bar{J}$ of $X$ commute (otherwise $J$ is not unique) and therefore $\bar{J}$ acts on the Riemann sphere $X/J$.
Once the reflection of the surface acts with fixed points, so does the induced reflection of the sphere. For a suitable choice of the global coordinate $x$ on the sphere, the reflection of the latter works as a complex conjugation.
It is convenient to think of $x$ as of the degree two meromorphic function on the surface. Obviously, it maps each real oval to the subset of the real equator $\hat{\mathbb{R}}:=\mathbb{R}\cup\infty$ of the Riemann sphere. There are three possibilities for the image of a real oval: (i) $x$ maps the oval $1-1$ to the equator; (ii) the oval is mapped $2-1$ to the equator or (iii) an oval is mapped $2-1$ to the finite segment of it. In the latter case the endpoints of the segment are the critical values of the projection $x$.

A simple combinatorial arguments show  that in the case of three real ovals
   on $X$, they are projected to three non-overlapping intervals of the real equator of the sphere. The endpoints of those intervals are the projections of the stationary points of the involution $J$ of the surface. We can give them the unique names in the following way.

Projections of the first, second and third oval induce either natural cyclic order of the real equator, or the inverse. In the latter case we change the sign of the coordinate $x$. Now we can enumerate the endpoints of the intervals in the natural cyclic order as $x_1,x_2,\dots,x_6$ so that  the first real oval of the surface is $x^{-1}([x_2,x_3])$, the second is $x^{-1}([x_4,x_5])$ and the third is $x^{-1}([x_6,x_1])$. Note that an interval of the extended real axis may contain infinity as its interior or boundary point. We have shown that our first definition of the moduli space is equivalent to the following

\begin{dfn} The moduli space  ${\cal M}_2\mathbb{R}_3$ is the factor  of   the cyclically ordered sextuples   of  points $(x_1, \dots, x_6)$ from
$\hat{\mathbb{R}}$ modulo the action of $PSL_2(\mathbb{R})$ (= real projective transformations conserving the orientation of the equator) .
\end{dfn}

 Normalizing the branch points e.g. as $x_4=\infty$, $x_5=-1$, $x_6=0$, one gets  the global  coordinate system on the moduli space:
\be
\label{moduli1}
0<x_1<x_2<x_3<\infty.
\ee
Other normalizations bring us to other global coordinate systems in the same space. Yet another global coordinate system in this space is related to the periods of holomorphic differentials.

\begin{thrm}
The period mapping $\Omega(X)$  is real analytic diffeomorphism from the moduli space ${\cal M}_2\mathbb{R}_3$  to the trihedral cone
\be
\label{PeriodsCone}
0<\Omega_{12}<min(\Omega_{11}, \Omega_{22})
\ee
\end{thrm}
{\bf Proof} (sketch).
First of all, the matrix $\Omega:=Im~\Pi$ is a well defined function on the
moduli space: both choices of the intrinsic basis in the space of integer 1-cycles introduced in Sect. \ref{SectBasisCycles} bring us to the same period matrix.

The CS integrals $u_j(x):=\int_{x_1}^x du_j$
map the upper half-plane $1-1$ to the rectangles with slots shown in the Fig. \ref{URect}. The dimensions of the rectangles are related to the elements of the period matrix where from the inequalities
(\ref{PeriodsCone}) follow.

Now we see that all the curves $X\in{\cal M}_2\mathbb{R}_3$ with fixed first column of the periods matrix
are parametrized by the length $l$ of the slot in the left rectangle of the Fig. \ref{URect}. In particular, $\Omega_{22}$
is a monotonic function of $l$. Now one can study the asymptotical behavior of this function and show that $\Omega_{22}(l)\to\infty$ when  $l\to 0$ \cite{Le2} and $\Omega_{22}(l)\to\Omega_{21}$ when $l\to1$ \cite{Le1}. ~~\bl

\begin{rmk}
The inverse  mapping, from the period matrices to the (suitably normalized) branchpoints of the curve is also effective. For genus two it is implemented by the Rosenhain formulae \cite{R} (see section \ref{theta}) in terms of theta constants.
\end{rmk}

\begin{rmk} The cone (\ref{PeriodsCone}) is strictly smaller than the space of all real positive definite symmetric $2\times2$ matrices $\Omega$. It is known (see e.g. \cite{Dub}) that any \emph{indecomposable} matrix from Siegel genus two space is a period matrix of some Riemann surface. The decomposable matrices make up a codimension one \emph{Humbert variety} which is determined by vanishing of at least one of ten even theta constants.  In our case of real curves the Humbert variety is a real codimension one subvariety in the real 3-space and its complement is disconnected. One can check that
\end{rmk}

\begin{tabular}{l}
$\theta\left[
{\tiny
\begin{array}{c} 11\\11\end{array}}
\right](i\Omega)=0$ on the edge $\{\Omega_{12}=0\}$ of the cone, \\
$\theta\left[
\tiny{
\begin{array}{c} 01\\10\end{array}}
\right](i\Omega)=0$ on the edge $\{\Omega_{12}=\Omega_{11}\}$ and\\
$\theta\left[
{\tiny
\begin{array}{c} 10\\01\end{array}}
\right](i\Omega)=0$ on the edge $\{\Omega_{12}=\Omega_{22}\}$.
\end{tabular}

The definitions of the theta constants will be given in the Sect. \ref{theta}.

\subsection{Curves with a marked point on real oval}
The problems of conformal mapping use a slightly more sophisticated moduli space, that of the genus two real curves with a marked point on a real oval (see e.g.  \cite{Bog1}).

\begin{dfn}
The space of genus two Riemann surfaces with three enumerated real ovals and a marked point $p_0\neq Jp_0$ on the third  oval we call ${\cal M}_{2,1}\mathbb{R}_3$. Two surfaces are equivalent iff there is a conformal mapping between them commuting with the reflections and respecting the enumeration of the ovals as well as the marked point.
\end{dfn}

An argument similar to that in the previous subsection shows that we can introduce an equivalent but more constructive

\begin{dfn}
By ${\cal M}_{2,1}\mathbb{R}_3$ we mean the sets of seven cyclically ordered  points $(x_0,x_1,\dots,x_6)$ in the real equator $\hat{\mathbb{R}}$ of the Riemann sphere modulo the action of $PSL_2({\mathbb R})$.
\end{dfn}
Here $x_0$ means the projection of the marked point to the sphere, other coordinates are the projections of the branchpoints of the curve.
The natural projection of ${\cal M}_{2,1}\mathbb{R}_3$ to the space ${\cal M}_2\mathbb{R}_3$ consists in forgetting of the marked point $x_0$.

Recall that for the element $X$ of the space ${\cal M}_2\mathbb{R}_3$
there is no natural distinction between two components of the set  $x^{-1}\mathbb{H}\subset  X$ until we orient a real oval.  The difference arises once we mark a point $p_0\neq Jp_0$ on a real oval: there is a unique disc $\mathbb{H}^+\subset$ $x^{-1}\mathbb{H}$ with $p_0$ on its boundary.

One can introduce several coordinate systems on the space ${\cal M}_{2,1}\mathbb{R}_3$. Fixing three points of seven, say $x_4:=\infty$,
$x_5:=-1$, $x_6:=0$, the positions of the remaining four points of the 7-tuple  will give us the global coordinate system on the moduli space:
\be
0<x_0<x_1<x_2<x_3<\infty.
\label{moduli2}
\ee
Other normalizations of the 7-tuple of points result in different coordinate systems on the moduli space. It is a good exercise to show that the arising coordinate change is a real analytic 1-1 mapping from the 4D cell (\ref{moduli2}) to the appropriate cell.

Yet another coordinate system on ${\cal M}_{2,1}\mathbb{R}_3$ is the modification of that related to the periods matrix. Three variables
$\Omega_{11}$, $\Omega_{12}$, $\Omega_{22}$ are inherited from the space
${\cal M}_{2}\mathbb{R}_3$ and the fourth is either $u_1^0$ or $u_2^0$,
the component of the image of the marked point under AJ map.
The integration path for $u(p_0)$ is the interval of the third real oval
avoiding the branch point $p_6$.

\begin{lmm}
Each mapping $(x_0,x_1,x_2,x_3) \to (\Omega_{11}, \Omega_{12}, \Omega_{22}, u_s^0)$, $s=1,2$ is a real analytic diffeomorphism of the cone (\ref{moduli2}) to the product of the cone (\ref{PeriodsCone}) and the interval $(0,1/2)$
\end{lmm}
{\bf Proof.}  Let us consider the points $x_s$, $s=0,1,2,3$ as the coordinates in the space ${\cal M}_{2,1}\mathbb{R}_3$, other branchpoints being fixed. By definition, the period matrix $\Omega$ is independent of the position of the marked point
$p_0$ and the mapping $(x_1,x_2,x_3)\to$ $\Omega$ is real analytic diffeomorphism from the cell (\ref{moduli1}) to the cone (\ref{PeriodsCone}). Both basic differentials $du_1$, $du_2$  are real and have no zeroes on the third real oval containing the marked point. Therefore, both values $u_s(p_0)$ monotonically increase from zero to
$\frac12=\frac12\int_{a_1+a_2}du_s$ when the marked point $p_0$ moves from $p_1$ to $p_6$ along the third real oval.
\bl

The periods of more sophisticated differentials taken instead of $du_s$ bring us to yet another coordinate systems in the same moduli space of curves with marked points on them.

\subsubsection{Christoffel-Schwarz differentials}\label{CSmap}
Let $1\le\alpha<\beta\le6$ be a couple of integers labeling the spaces of heptagons. To each element of the moduli space
${\cal M}_{2,1}\mathbb{R}_3$ we ascribe the unique abelian differential $dw_{\alpha\beta}$ of the third kind with simple poles at the points $p_0, Jp_0$, residues resp. $-1,+1$ at those points and zeros at the branchpoints $p_\alpha$ and $p_\beta$ (one of the ways to normalize meromorphic differentials is to ascribe its zeroes at the points of nonspecial degree $g$ divisor). The differential will automatically have double zeroes because is is odd with respect to the involution $J$. In the algebraic model (\ref{X}) it takes the form:
\be
\label{dwab}
dw_{\alpha\beta}=(x-x_\alpha)(x-x_\beta)\frac{dx}y
\ee
where $y^2=\prod_{j=1}^6(x-x_j)$ and the 7-tuple $\infty=:x_0,x_1,\dots,x_6$
represents an element of the space ${\cal M}_{2,1}\mathbb{R}_3$.
Christoffel-Schwarz differential (\ref{dwab}) can be decomposed into a sum of three elementary ones:
\be
dw_{\alpha\beta}=dv_{Jp_0p_0}+C_1du_1+C_2du_2
\label{diffdecomp}
\ee
where $dv_{Jp_0p_0}$ is $a$-normalized third kind abelian differential with poles at $Jp_0$ and $p_0$; $du_1$ and $du_2$ are above normalized holomorphic differentials and the constants $C_1$ and $C_2$ are uniquely determined by the condition that $dw_{\alpha\beta}$ has zeroes at the branching points $p_\alpha$ and $p_\beta$.

Each element of the heptagon space ${\cal P}_{\alpha\beta}$ is the image of the upper half plane under the Christoffel-Schwarz map
\be
\label{CS}
w_{\alpha\beta}(x):=\int_*^xdw_{\alpha\beta}.
\ee
Actually, even more is true.

\begin{figure}[h]
\begin{picture}(130,65)
\thicklines
\put(60,60){\line(1,0){60}}
\put(60,60){\vector(0,-1){35}}
\put(60,25){\line(-1,0){55}}
\put(5,25){\line(0,-1){25}}
\put(5,0){\line(1,0){40}}
\put(45,0){\vector(0,1){45}}
\put(45,45){\line(1,0){75}}
\put(52,35){$*$}
\end{picture}
\caption{\small Self-intersection of the image of
 real equator under $w_{26}$ map.}
\label{SelfIS}
\end{figure}
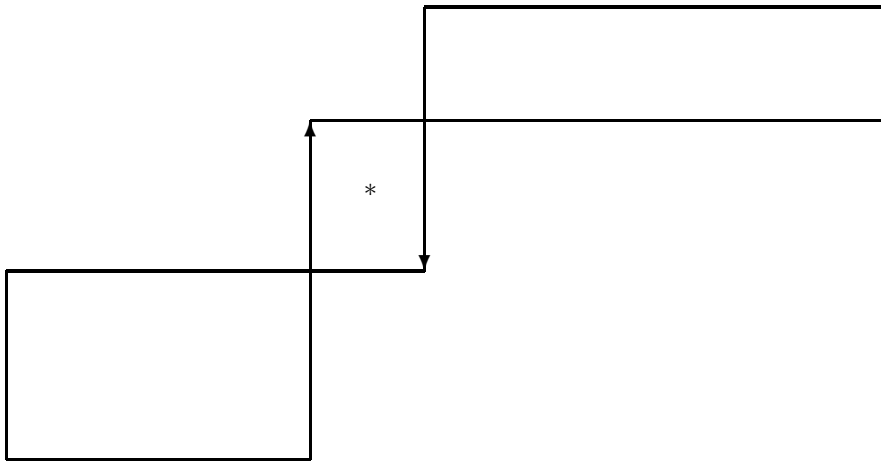

\begin{thrm}
\label{th2}
Christoffel-Schwarz mapping $w_{\alpha\beta}(x)$
induces a real analytic diffeomorphic map from the moduli space
${\cal M}_{2,1}\mathbb{R}_3$ to the heptagon space ${\cal P}_{\alpha\beta}$.
\end{thrm}
{\bf Proof.}
First of all we check that each CS differential $dw_{\alpha\beta}$ is real.
This implies that the increment of the CS map on the boundary of the disc $\mathbb{H}^+$ embedded to the Riemann surface is real on the real ovals  and it is pure imaginary on the coreal ovals. This increment is monotonic
between the branchpoints but at $p_0$. The image of $\partial\mathbb{H}^+$ under the CS map is a polygonal line with the same sequence of corners as in any heptagon from the space ${\cal P}_{\alpha\beta}$, namely all the corners are equal to $\pi/2$ but those at $p_\alpha$ and $p_\beta$ where it is equal to $3\pi/2$. Moreover, this polygonal line has no self-intersections.
For instance, in the case $(\alpha,\beta)=(2,6)$ the self-intersection shown in  the Fig. \ref{SelfIS} is impossible because the points of  the rectangle marked by  $*$ have $-1$ preimages in $\mathbb H$ according to the argument principle.

We have established  a $1-1$ correspondence between the moduli space and the space of heptagons induced by the CS map (\ref{CS}). Now we prove the real analyticity of the direct and the inverse mapping. One can easily see  that  the dimensions of the heptagon are the periods of the corresponding CS integral (since CS differential is odd with respect to $J$):
\be
\label{auxparam1}
\begin{array}{l}
iH_1:=\int_{p_2}^{p_1}dw_{\alpha\beta}=\frac12\int_{b_1}dw_{\alpha\beta};\\
i^2H_2:=\int_{p_3}^{p_2}dw_{\alpha\beta}=-\frac12\int_{a_1}dw_{\alpha\beta};\\
i^3H_3:=\int_{p_4}^{p_3}dw_{\alpha\beta}=i\pi+\frac12\int_{b_2-b_1}dw_{\alpha\beta};\\
i^4H_4:=\int_{p_5}^{p_4}dw_{\alpha\beta}=-\frac12\int_{a_2}dw_{\alpha\beta};\\
i^5H_5:=\int_{p_6}^{p_5}dw_{\alpha\beta}=\frac12\int_{-b_2}dw_{\alpha\beta}.
\end{array}
\ee
The basic 1-cycles may be separated from the branchpoints of the surface, so the real analyticity of the direct map ${\cal M}_{2,1}\mathbb{R}_3\to $  ${\cal P}_{\alpha\beta}$ is clear. It remains to show that the map has the full rank.

Let us consider the coordinate system in the moduli space such that
$x_0=\infty$, $x_\alpha=-1$, $x_\beta=1$. Assume that  the CS-induced map degenerates at a point $(x_0,x_1,\dots,x_6)$ of the moduli space, then there exist a nontrivial tangent vector $\xi:=\sum_{1=j\not\in\{\alpha,\beta\}}^6
\xi_j\frac\partial{\partial x_j}$ annihilating all the periods of the CS differential $dw_{\alpha\beta}$
at this point. This means the existence of the meromorphic differential
$$
dv:= \frac12\sum_{1=j\not\in\{\alpha,\beta\}}^6
\xi_j\frac{dw_{\alpha\beta}}{x-x_j}
$$
with zero cyclic and polar periods on the associated surface and with poles at the branchpoints of the surface, except $p_\alpha$ and $p_\beta$.
Then the meromorphic function $v(x,y):=\int_{p_\alpha}^{(x,y)}dv$ is single valued on the surface and it has at most four simple poles along with triple zeroes at the points $p_\alpha$ and $p_\beta$ (since $\int_\alpha^\beta dv=0$). This function is identical zero and therefore the above tangent vector is also zero. Hence the CS-induces mapping ${\cal M}_{2,1}\mathbb{R}_3\to $  ${\cal P}_{\alpha\beta}$ has the full rank everywhere. ~~~\bl

\subsubsection{Auxiliary parameters of the CS map}
Given a heptagon from the space ${\cal P}_{\alpha\beta}$, the corresponding point in the moduli space ${\cal M}_{2,1}\mathbb{R}_3$ may be found as (unique as it follows from the Theorem \ref{th2}) solution of the following system of equations obtained from (\ref{auxparam1}) by substituting the decomposition (\ref{diffdecomp}) of Christoffel-Schwarz differential into elementary differentials. (We have also used
Riemann bilinear relations (\ref{RB3}) to transform the $b$-periods of elementary differential of the 3-rd kind to Abel-Jacobi image of its poles.)
\be
\label{auxparam2}
\begin{array}{l}
2H_1=(C_1\Omega_{11}+C_2\Omega_{12}+2\pi(1-2u_1^0));\\
2H_2=C_1;\\
2H_4=-C_2;\\
2H_5=-(C_1\Omega_{12}+C_2\Omega_{22}-4\pi u_2^0).
\end{array}
\ee
Here $C_1$ and $C_2$ are real analytic functions on the moduli space
defined above as a solution of linear
system with the nonsingular $2\times2$ matrix $||du(p_\alpha),du(p_\beta)||$.
We see that the auxiliary system of equations is "almost linear"
with respect to the coordinate system related to the periods of holomorphic differentials. It remains to get the effective evaluation of the functions
$C_1,C_2$ and $u_2^0$ as well as the CS integral itself.
This is the subject of the next section.

\section{Theta functions on genus two surfaces}\label{theta}
Here we give a short introduction to the theory of Riemann theta functions adapted to genus two surfaces. Three problems related to conformal mappings will be effectively solved in terms of Riemann theta functions:
\begin{itemize}
\item Localization of the curve inside its Jacobian;
\item Representation of the 2-sheeted projection of the curve to the sphere;
\item Evaluation of the normalized abelian integral of the third kind (which is the essential part of the CS-integral).
\end{itemize}

\begin{dfn}
Let $u\in\mathbb{C}^2$ and $\Pi\in\mathbb{C}^{2\times2}$ be a Riemann
matrix, i.e. $\Pi=\Pi^t$ and $Im~\Pi>0$. The theta function of those two arguments is the following Fourier series

\be
\label{thetadef}
\theta(u, \Pi):=\sum\limits_{m\in\mathbb{Z}^2}
\exp(2\pi im^tu+\pi i m^t\Pi m),
\ee
\end{dfn}

The convergency of this series grounds on the positive determinacy of
$Im~\Pi$.  The series has high convergency rate with well controlled accuracy \cite{DHB}.
Theta function has the following easily checked quasi-periodicity properties with respect to the lattice $L(\Pi):=\Pi\mathbb{Z}^2+\mathbb{Z}^2$:
\be
\label{quasiperiod}
\theta(u+\Pi m+m', \Pi)=\exp(-i\pi m^t\Pi m-2i\pi m^tu)\theta(u,\Pi),
\qquad m,m'\in\mathbb{Z}^2.
\ee

   One can say that theta function is a section of a certain line bundle over the Jacobian. In particular, its zeroes make up a well defined set known as a theta divisor $(\theta)$.
The Abel-Jacobi map transfers the theta function to the Riemann surface where it becomes the multivalued function
(or a section of a certain line bundle)
$$
\theta_e(p):=\theta(u(p)-e;\Pi),
$$
which is multiplied by some nonvanishing factors once its argument $p$ goes around the handles of the surface. Therefore, its zeroes are well defined.

The zero set of the theta function is described by so called Riemann vanishing theorems the most important of them is the following \cite{FK}.

\begin{thrm}[Riemann]
Function $\theta_e(p)$ on a genus two surface either\\
(i) vanishes identically on the surface iff $e$
is a certain effectively calculated point  $K$ of Jacobian
(a.k.a. vector of Riemann constants), or\\
(ii) has exactly two zeroes $q_1,q_2$  such that
$q_1\neq Jq_2$ and
\be
\label{RvT}
e=u(q_1)+u(q_2)+K ~~mod~L(\Pi)
\ee
\end{thrm}

\subsection{Vector of Riemann constants}
We have to determine the vector of Riemann constants $K$ for our particular choice of the initial point in the AJ map and the choice of the basis in the homologies. Consider the theta function with characteristics
which is the slight modification of the above theta.
\begin{dfn}
\be
\label{thetachardef}
\theta[2\epsilon, 2\epsilon'](u, \Pi):=\sum\limits_{m\in\mathbb{Z}^2}
\exp(2\pi i(m+\epsilon)^t(u+\epsilon')+\pi i (m+\epsilon)^t\Pi (m+\epsilon))
\ee
$$
=
\exp(i\pi\epsilon^t\Pi\epsilon+2i\pi\epsilon^t(u+\epsilon'))
\theta(u+\Pi\epsilon+\epsilon',\Pi),
\qquad \epsilon,\epsilon'\in\mathbb{R}^2.
$$
The matrix argument $\Pi$ of theta function is usually omitted when it is clear which matrix we mean.  Omitted vector argument $u$ is supposed to be zero and the appropriate function of $\Pi$ is called the theta constant:
$$
\theta[\epsilon, \epsilon']:=\theta[\epsilon, \epsilon'](0, \Pi).
$$
\end{dfn}

\begin{rmk}\label{RemTheta}
(i) Theta function with integer characteristics
$[2\epsilon, 2\epsilon']$ is either even or odd depending on the parity of the inner product $4\epsilon^t\cdot\epsilon'$. In particular, all odd theta constants are zeroes.\\
(ii) Adding matrix with even entries to integer theta characteristics can
at most spoil the sign of the theta function. Hence, the binary arithmetic plays a great role in the calculus of theta functions.
\end{rmk}

\begin{lmm}
\be
\label{K}
K=u(p_3)+u(p_5)
\ee
\end{lmm}
{\bf Proof.}
We see from the table is subsection \ref{AJPj} that the points $p_3$ and $p_5$ correspond to the odd characteristics. Remark \ref{RemTheta}(i) implies that the function $\theta(u(p))$ will have zeroes at the points $p=p_3,p_5$ on the surface.
This function cannot be identical zero on the surface:
otherwise from part (i) of Riemann vanishing theorem it would follow that $K=0$. Now one can check that the function $\theta_e(p)$
with the shift $e:=u(p_3)+u(p_5)+K$ is not identical zero (since $p_3\neq Jp_5$) and vanishes at three points $p=p_1,p_3,p_5$.

An alternative (ii) of Riemann vanishing theorem suggests that $K=u(p_3)+u(p_5)$ $mod ~L(\Pi)$. In other words, $K$ corresponds to the odd characteristics
$\tiny \left[\begin{array}{c} 10\\11\end{array}\right]$. ~~~\bl

\begin{rmk}
It is convenient to represent integer theta characteristics as the sums of AJ images of the branchpoints, keeping only the indexes of those points:
$[sk..l]$ means the sum modulo 2 of the theta characteristics of points
$u(p_s)$, $u(p_k)$, $\dots,u(p_l)$, e.g. $[35]$ is the vector of Riemann constants represented by theta characteristics.
\end{rmk}

\subsection{Image of Abel-Jacobi map}
The location of genus 2 curve embedded to its Jacobian may be reconstructed by solving a single equation.

\begin{thrm}[Riemann]
A point $e$ of Jacobian lies in the image $u(X)$ of Abel-Jacobi map
if and only if
$$
\theta{\tiny\left[\begin{array}{c}10\\11\end{array}\right]}(e)=0
$$.
\end{thrm}
{\bf Proof.} 1.  The function $\theta(u(p)-K)$ vanishes identically on the curve and so does the function $\theta[35](u(p))$.

2. Conversely, suppose that $\theta[35](e)=0$ or equivalently, $\theta(e+K)=0$.
The function $\theta_{e'}(p)$ with the shift $e':=e+K$
either vanishes at two points $p_1$, $p'$ and in this case $e=u(p')$,
or identically, then $e=0=u(p_1)$. ~~~\bl

\subsection{Projection to the sphere}\label{projection}
Any meromorphic function on the curve my be effectively calculated via the
Riemann theta functions once we know its divisor. Take for instance the degree 2 function $x$ on the hyperelliptic curve (\ref{X}). This projection is unique if normalized  e.g. as follows:
$x(p_s)=0$, $x(p_j)=1$, $x(p_l)=\infty$, $s,j,l$ is a positive triple from the index set $\{1,2,\dots,6\}$.

With the use of the transformation rules (\ref{quasiperiod}) one immediately checks that the following function is single valued on the curve:
$$
\tilde{x}(p):= \frac{\theta^2[sk35](u(p),\Pi)}{\theta^2[lk35](u(p),\Pi)},
\qquad k\neq s,l.
$$
Now with the help of Riemann theorem
one checks that the numerator of the function has double zeros at the points $p_s,~p_k$ while the denominator has double zeroes at $p_l,~p_k$.
Therefore, $\tilde{x}(p)$ differs by a constant factor from the above normalized projection:
\be
\label{xofP}
x(p)= \pm\frac{\theta^2[lkj35]}{\theta^2[skj35]}
      \frac{\theta^2[sk35](u(p))}{\theta^2[lk35](u(p))},
\qquad k\neq s,l,j,
\ee
where the sign $\pm$ in the latter formula
depends on the parity ($+$ even /$-$ odd)
of the scalar product $\epsilon(j)\cdot(\epsilon'(s)+\epsilon'(l))$, where
$\Pi\epsilon(s)+\epsilon'(s):=u(p_s)$. In the case $k=j$ the normalization have to use the L'Hospital rule and we get a slightly different factor.

Also, the standard normalization $x(p_\alpha)=-1; x(p_\beta)=1; x(p_0)=\infty$
brings us to a little more awkward expression:
$$
x(p):=\pm2\frac{\theta^2[\beta j35](u^0)}{\theta^2[\alpha\beta j35]}
\frac{\theta^2[\alpha j35](u(p))}{\prod_\pm\theta[j35](u(p)\pm u^0)}-1,
\qquad j\neq\alpha,\beta,
$$
again, the sign in front of the fraction depends on the parity
of the scalar product $\epsilon(\beta)\cdot\epsilon'(\alpha)$, where
$\epsilon(s)$, $\epsilon'(s)$ is the representation of the half-period
$u(p_s)$ by integer theta characteristic.

\subsection{Rosenhain formulae; Humbert variety}
Putting in the above formula (\ref{xofP}) the half-periods corresponding to the branchpoints, we get the effectively computed expressions for the latter in terms of the period matrix. For instance, if we normalize the projection $x(p)$
as $x(p_1)=0$, $x(p_2)=1$, $x(p_6)=\infty$, we get the following expressions for the remaining branch points in terms of theta constants:
\be
x_3=\frac{
\theta^2\left[
{\tiny
\begin{array}{c} 00\\00\end{array}}
\right]
\theta^2\left[
{\tiny
\begin{array}{c} 00\\01\end{array}}
\right]
}
{
\theta^2\left[
{\tiny
\begin{array}{c} 01\\00\end{array}}
\right]
\theta^2\left[
{\tiny
\begin{array}{c} 01\\01\end{array}}
\right]
};
\quad
x_4=\frac{
\theta^2\left[
{\tiny
\begin{array}{c} 00\\01\end{array}}
\right]
\theta^2\left[
{\tiny
\begin{array}{c} 10\\10\end{array}}
\right]
}
{
\theta^2\left[
{\tiny
\begin{array}{c} 01\\00\end{array}}
\right]
\theta^2\left[
{\tiny
\begin{array}{c} 11\\11\end{array}}
\right]
};
\quad
x_5=\frac{
\theta^2\left[
{\tiny
\begin{array}{c} 00\\00\end{array}}
\right]
\theta^2\left[
{\tiny
\begin{array}{c} 10\\10\end{array}}
\right]
}
{
\theta^2\left[
{\tiny
\begin{array}{c} 11\\11\end{array}}
\right]
\theta^2\left[
{\tiny
\begin{array}{c} 01\\01\end{array}}
\right]
}.
\label{branch1}
\ee
Choosing other normalization for the projection, we get certain expressions for the cross-ratios of the same set of branch points, e.g.
\be
1-x_3=-\frac{
\theta^2\left[
{\tiny
\begin{array}{c} 10\\00\end{array}}
\right]
\theta^2\left[
{\tiny
\begin{array}{c} 10\\01\end{array}}
\right]
}
{
\theta^2\left[
{\tiny
\begin{array}{c} 01\\00\end{array}}
\right]
\theta^2\left[
{\tiny
\begin{array}{c} 01\\01\end{array}}
\right]
};
\quad
1-x_4=-\frac{
\theta^2\left[
{\tiny
\begin{array}{c} 00\\10\end{array}}
\right]
\theta^2\left[
{\tiny
\begin{array}{c} 10\\01\end{array}}
\right]
}
{
\theta^2\left[
{\tiny
\begin{array}{c} 01\\00\end{array}}
\right]
\theta^2\left[
{\tiny
\begin{array}{c} 11\\11\end{array}}
\right]
};
\quad
1-x_5=-\frac{
\theta^2\left[
{\tiny
\begin{array}{c} 00\\10\end{array}}
\right]
\theta^2\left[
{\tiny
\begin{array}{c} 10\\00\end{array}}
\right]
}
{
\theta^2\left[
{\tiny
\begin{array}{c} 11\\11\end{array}}
\right]
\theta^2\left[
{\tiny
\begin{array}{c} 01\\01\end{array}}
\right].
}
\label{branch2}
\ee
Comparing the formulae for the same branch points from the equations (\ref{branch1}) and (\ref{branch2}) we get certain relations for the theta constants which  are the consequences from the Riemann theta identities.
The expressions for the branch points of genus two (hence hyperelliptic) curves in terms of theta constants appeared in Rosenhain's work \cite{R}.

\begin{thrm} \cite{R}
For nonsingular genus 2 curve all 10 even theta constants
$\theta[\epsilon,\epsilon']$, $\epsilon^t\cdot\epsilon'\in2\mathbb{Z}$, are distinct from zero.
\end{thrm}
{\bf Proof.} Choose three distinct branching points $p_s$, $p_j$, $p_l$ on the curve.
Since $p_s\neq Jp_j$, the function  $\theta_e(p)$ does not vanish at $p=p_l$
when $e=u(p_s)+u(p_j)+K$. The value $\theta_e(p_l)$ vanishes simultaneously  with the theta constant $\theta[sjl35]$. Once the indexes $s,j,l$ vary in the set $\{1,2,\dots,6\}$, the theta characteristic
$[sjl35]$ runs through all 10 even characteristics. ~~\bl

Humbert variety in complex Siegel space is the locus where at least one of even theta constants vanishes. This is exactly the locus of Riemann matrices
that are \emph{not} period matrices of any genus two Riemann surface.

\subsection{Third kind abelian itegral}
On any Riemann surface there exist a unique abelian differential  of the third kind $dv_{rq}$ with simple poles at two prescribed points $r,q$ only, residues $+1,-1$ respectively and trivial periods along all $a-$ cycles. It has certain physical meaning in terms of flow of inviscid incompressible fluid on the surface with a source at $r$ and sink at $q$. Our approach to effective conformal mapping grounds on the fact that the integral of this differential
can be expressed in closed form.

\begin{thrm}[Riemann]
Choose any point $s\in X$, then for any two points $r,q\neq Js$ of the surface holds the representation:
\be
v_{rq}(p):=\int_*^pdv_{rq}=
\log\frac{\theta(e+u(p)-u(r))}{\theta(e+u(p)-u(q))}=
\log\frac{\theta[\epsilon,\epsilon'](u(p)-u(r))}
{\theta[\epsilon,\epsilon'](u(p)-u(q))}+const,
\label{intrep}
\ee
where theta characteristic
$[\epsilon,\epsilon']$ corresponds to the zero
$e:=-u(s)-K$ of theta function (say, it may be any odd integer one).
\end{thrm}
{\bf Proof.}
Consider the following function of variable $p\in X$:
$$
\exp(-v_{rq}(p))
\frac{\theta(e+u(p)-u(r))}{\theta(e+u(p)-u(q))}.
$$
One checks that it is locally holomorphic on the surface: the theta function in numerator vanishes at the points $r,s$, the denominator vanishes at $q,s$.
Moreover, this function is single valued on the surface\footnote{provided the branches for $u(r),u(q)$ in (\ref{intrep}) are properly chosen:
the integration path from $r$ to the initial point $p_1$ and further to $q$ may be deformed to the one disjoint from $a-$ and $b-$ cycles}. To prove this we use the Riemann bilinear identity
\be
\label{RB3}
\int_{b_j}dv_{rq}=2\pi i\int_q^rdu_j
\ee
and the transformation properties (\ref{quasiperiod}) of theta functions.
Therefore, it is a constant independent of $p$. ~~\bl

{\bf Example:} CS integral from the section \ref{CSmap} may be represented as
follows:
$$
w_{\alpha\beta}(p)=log\frac{\theta[3](u(p)+u^0)}{\theta[3](u(p)-u^0)}
+C_1u_1(p)+C_2u_2(p),
$$
where $u(p):=(u_1(p),u_2(p))^t$; $u^0:=u(p_0)$ and the constants $C_1,C_2$ are obtained from the system of linear equations
$dw_{\alpha\beta}(p_\gamma)=0$, when $\gamma=\alpha,~\beta$.

\section{Algorithm of conformal mapping}
Based on the formulae of the previous section, we can propose the algorithm of conformal mapping of the heptagon to the half-plane and vice versa. First of all, given the heptagon we have to determine the corresponding point of the moduli space ${\cal M}_{2,1}\mathbb{R}_3$.

\subsection{Auxiliary parameters}
Given (related to side lengths) coordinates $H_s$ of the heptagon, we have to determine seven real parameters: the (imaginary part of period) matrix $\Omega$, the image $u^0:=u(p_0)$ of the marked point $p_0$ in the Jacobian of the curve and real 2-vector $C:=(C_1,C_2)$. Those parameters give a solution to a system of seven real equations:

\be
\label{dweq0}
d\theta[35](u,i\Omega)\wedge
d(\log\frac{\theta[3](u+u^0,i\Omega)}{\theta[3](u-u^0,i\Omega)}+C\cdot u)=0, \qquad {\rm when}~u=u(p_\alpha), u(p_\beta)
\ee
which means that CS differential $dw_{\alpha,\beta}$ has zeroes at the points $p_\alpha$ and $p_\beta$;
\be
\label{ovalu0}
\theta[35](u^0,i\Omega)=0,
\ee
which means that the point $u^0$ lies on the AJ image of the curve in the Jacobian and finally
\be\label{Sides}
\begin{array}{l}
H_1=\frac12(C_1\Omega_{11}+C_2\Omega_{12}+2\pi(1-2u_1^0));\\
H_2=\frac12  C_1;\\
H_4=-\frac12 C_2;\\
H_5=-\frac12(C_1\Omega_{12}+C_2\Omega_{22}-4\pi u_2^0).
\end{array}
\ee
which specify the side lengths of the heptagon.

\begin{lmm}
Let the side lengths $H_1,H_2\dots,H_5$ satisfy the restrictions described in  lemma \ref{H1H5}, then the system of seven equations (\ref{dweq0}), (\ref{ovalu0}), (\ref{Sides}) has a unique solution $\Omega$, $u^0$, $C\in\mathbb{R}^2$ in the domain determined by inequalities:
$$
0<\Omega_{12}<\min(\Omega_{11},\Omega_{22}),
$$
$$
0<u_s^0<1/2, \qquad s=1,2.
$$
\end{lmm}
{\bf Proof.} The existence and the uniqueness of the solution to (\ref{dweq0}) -- (\ref{Sides}) in the specified domain follows from the existence and the uniqueness (up to real projective transformations) of the conformal mapping of a given heptagon to the upper half plane. ~~~\bl

\begin{rmk}
Essentially we have just three nonlinear equations to solve which can be done by Newton method with parametric continuation.
Indeed, take any point $(\Omega$, $u_1^0)$ in the moduli space, then solving just one nonlinear equation (\ref{ovalu0}) we get the whole vector $u^0$, then solving two linear equations  (\ref{dweq0}) we get the constants $C$ and  substituting the values $\Omega$, $u^0$, $C$ to the last four equations (\ref{Sides}) we get the dimensions of the appropriate heptagon. So we got the correspondence of two points: one in the space ${\cal P}_{\alpha,\beta}$
and another -- in the space ${\cal M}_{2,1}\mathbb{R}_3$. Now, any given point in the heptagon space may be connected to the obtained one by a path of at most two linear segments. Those may be effectively lifted to the moduli space by the Newton method.
\end{rmk}

\subsection{Mapping heptagon to half plane}
 Suppose a point $w^*$ lies in the normalized heptagon (i.e. with the vertex  $w_1=i\pi$). Consider a system of two equations
\be
\label{heptoH}
\begin{array}{r}
\log\frac{\theta[3](u^*+u^0,i\Omega)}{\theta[3](u^*-u^0,i\Omega)}+C\cdot u^* = w^*,\\
\theta[35](u^*,i\Omega)=0,
\end{array}
\ee
with respect to the unknown complex 2-vector $u^*$;
real symmetric matrix $\Omega$ and real 2-vectors $u^0$, $C$ being the solutions of the auxiliary system (\ref{dweq0}) -- (\ref{Sides}).

 \begin{lmm}
  System (\ref{heptoH}) has the unique solution $u^*$ with theta characteristic from the block
  $\tiny\left[
\begin{array}{cc}
-I~&I\\
-I~&I\\
\end{array}
\right]$.
 \end{lmm}
{\bf Proof.} 
The left hand side of the first equation in (\ref{heptoH}) restricted to 
the AJ image of the marked surface $(X,p_0)$ in $\mathbb{C}^2$
represents the abelian integral $w_{\alpha\beta}$.
The abelian integral is not a single valued function on the surface $X$. However there is a single valued branch of the integral in the disc $\mathbb{H}^+$ which maps the latter $1-1$ to the given normalized heptagon. Therefore, the point $w^*$ has the unique preimage in the disc $\mathbb{H}^+\subset X$. Abel-Jacobi image $u^*$ of this point satisfies two equations (\ref{heptoH}) and lies in the specified block of the space $\mathbb{C}^2$. ~~~\bl

\begin{rmk}
It is clear from the reflection principle for the conformal mappings that the set of two equations (\ref{heptoH}) may have many
solutions. We use theta characteristic to single out the unique one.
\end{rmk}

Substituting the solution $u^*$ to the right hand side of the expression (\ref{xofP}) we get the evaluation at the point $w^*$ of the conformal mapping $x(w)$ of the heptagon to the half plane with normalization $w_s$, $w_j$, $w_l\to$  $0,1,\infty$.
\be
\label{xofP2}
x^*= \frac{\theta^2[lkj35]}{\theta^2[skj35]}
      \frac{\theta^2[sk35](u^*)}{\theta^2[lk35](u^*)},
\qquad k\neq s,l,j.
\ee

\subsection{Mapping half plane to heptagon}
Conversely, given a point $x^*$ in the upper half plane we solve a system of two equations with respect to a complex 2-vector $u^*$ with characteristics from
$\tiny\left[
\begin{array}{cc}
-I~&I\\
-I~&I\\
\end{array}
\right]$:
\be
\label{Htohept}
\begin{array}{r}
\frac{\theta^2[lkj35]}{\theta^2[skj35]}
      \frac{\theta^2[sk35](u^*)}{\theta^2[lk35](u^*)}=x^*,\\
\theta[35](u^*)=0
\end{array}
\ee
and substitute this solution to the formula
\be
\log\frac{\theta[3](u^*+u^0)}{\theta[3](u^*-u^0)}+C\cdot u^*:=w^*
\ee
to get the image of the point $x^*$ in the heptagon.

\section{Conclusion}
It often happens in classical mechanics that the problem becomes much simpler when considered in adequate coordinate system. We see that the same is true for the conformal mapping: taking elements of the period matrix as independent coordinates on the moduli space of curves, we have simplified the system of equations for the auxiliary parameters of conformal mapping.
The mappings themselves are now closed analytic  expressions involving effectively computed Riemann theta functions.

 Our approach to conformal mapping of (rectangular) polygons can readily be generalized to other types of polygons. Of course, for involved polygons the the genus of the associated curve becomes higher which makes the algorithm more complicated. For the set of auxiliary parameters one has to solve the (hyperelliptic) Schottky problem additionally. As well, more efforts have to be spent to localize the AJ image of the curve inside the multidimensional Jacobian.

\section{Appendix: asymptotic of periods}

\vspace{5mm}
\parbox{9cm}
{\it
119991 Russia, Moscow GSP-1, ul. Gubkina 8,\\
Institute for Numerical Mathematics,\\
Russian Academy of Sciences\\[3mm]
{\tt gourmet@inm.ras.ru, ab.bogatyrev@gmail.com}}


\begin{thebibliography}{15}
\bibitem{TD} L.N. Trefethen, T.A. Driscoll, Schwarz-Christoffel Mapping --      Cambridge Univ. Press, 2002.

\bibitem{Nat}  S.M.Natanzon,  Moduli of surfaces, real algebraic curves and and their superanalogs -- AMS Translation of Math. Monographs, 2004

\bibitem{RF} Rauch, H.E and H.M. Farkas,  Theta functions with applications to Riemann surfaces -- Williams \& Wilkins Company, Baltimore, 1974

\bibitem{FK} H.M. Farkas and I.Kra,  Riemann Surfaces
-- Springer Verlag, NY, Heidelberg, Berlin, 1980

\bibitem{DHB}
Deconinck,~B.; Heil,~M.; Bobenko,~A.; van Hoeij,~M.; and Schmies,~M. Computing Riemann Theta Functions// Math. Comput. 73, 1417-1442, 2004.

\bibitem{BHY} A.Bogatyrev, M.Hassner, D.Yarmolich  An exact analytical-expression for the read sensor signal in magnetic data storage channels// in ``Error-Correcting Codes, Finite Geometries and Cryptography'', eds. A.A.Bruen, D.L.Wehlau, AMS series Contemporary Math. 523 (2010),  155-160.

\bibitem{Dub} B.A.Dubrovin, Riemann Surfaces and Nonlinear Equations -- AMS, 2002

\bibitem{Bog1}  A.B.Bogatyrev, Effective approach to least deviation problems// Sbornik: Math, 2002, 193:12, 1749--1769.

\bibitem{Le2} A.Lebowitz, Handle removal on a compact Riemann surface of genus 2//
Israel Journal of Mathematics, 15:2, 189-192


\bibitem{Le1} A.Lebowitz, Degeneration of a compact Riemann surface of genus 2
// Israel Journal of Mathematics, 12:3, 223-236,

\bibitem{R} G.Rosenhain, Abhandlung \"uber die Funktionen zweier Variabeln mit Vier perioden, Mem. pres. l'Acad. de Sci. de France des savants XI (1851).

\end{thebibliography}
\end{document}